\documentclass[a4paper, 12pt]{article}
\usepackage{amsmath}
\usepackage{amsfonts}
\usepackage{amsthm}
\usepackage{epsfig}
\usepackage[pdftex]{hyperref}
\usepackage{indentfirst}

\newtheorem*{teo}{Theorem}
\newtheorem{thm}{Theorem}

\newtheorem{prop}{Proposition}
\newtheorem{cor}{Corollary}

\newcommand{\re}{\mathbb{R}}
\newcommand{\hp}{\mathbb{H}^{2}}
\newcommand{\hs}{\mathbb{H}^{3}}

\newcommand{\sw}{\mathbb{S}^{2}}
\newcommand{\st}{\mathbb{S}^{3}}
\newcommand{\si}{\Sigma}

\newcommand{\ric}{\operatorname{Ric}}

\newcommand{\dist}{\operatorname{dist}}

\title{Area Estimates and Rigidity of Non-compact $H$-Surfaces in 3-Manifolds}
\author{Vanderson Lima}
\date{}

\begin{document}

\maketitle

\begin{abstract}
For appropriately values of $H$, we obtain an area estimate for a complete non-compact $H$-surface of finite topology and finite area, embedded in a three-manifold of negative curvature. Moreover, in the case of equality and under additional assumptions, we prove that a neighbourhood of the mean convex side of the surface must be isometric to a hyperbolic Fuchsian manifold. Also, we show by an counter-example that although that area estimate holds for minimal surfaces, one does not have rigidity for equality in this case.

\end{abstract}

\providecommand{\abs}[1]{\lvert#1\rvert}

\linespread{1} % Espaçamento entre as linhas

\section{Introduction}

Let $\si$ be a closed $H$-{\it surface} (a surface of constant mean curvature $H$) embedded in a Riemannian $3$-manifold $(M,g)$. Imposing a curvature condition on $(M,g)$ and perhaps some extra condition on $\si$ is possible to obtain area estimates and rigidity results for $\si$ and $(M,g)$. For example, assuming a lower bound on the scalar curvature and that the surface is minimizing or does have index $1$ with respect to the functional $Area - H\cdot Volume$, some rigidity results were obtained in \cite{C.G, Y, B.B.E.N, B.B.N, N, M.N}. There are also related results for {\it free boundary minimal surfaces} in $3$-manifolds with boundary, see \cite{Am,Me1}, and for {\it MOTS} in a spacetime, see \cite{Ga.Me,Me2}. See also \cite{B.C,Mo} for results in higher dimensions.

The case of interest on this work is when one have bounds on the sectional curvature $K_{M}$ of $(M,g)$ but no assumption on the index of the $H$-surface $\si$. We will state two important results on this direction.
\begin{teo}[Mazet-Rosenberg, \cite{Ma.Ro}]
Let $(M^3,g)$ be a complete oriented Riemannian $3$-manifold whose sectional curvatures satisfy $0 \leq K_{M} \leq 1$. Let $\si \subset M$ be an embedded minimal $2$-sphere, then $$|\si| \geq 4\pi.$$ 
Moreover the equality holds if, and only if, $(M,g)$ is isometric to the round sphere $\st_{1}$ or to a quotient of $\sw_{1}\times\re$ with the round product metric.
\end{teo}

To state the second result consider the following construction.\\\\
\noindent
{\bf Fuchsian $3$-Manifolds}: Let $\si$ be a oriented surface which is either closed or non-compact with finite topology, endowed with a complete metric $h_{\si}$ of constant curvature $-1$ and finite area, so we have that $\chi(\si) < 0$. The Riemannian manifold $(\si,h_{\si})$ is a quotient of $\hp$ by a cocompact or cofinite Fuchsian group $\Gamma$ of isometries of $\hp$. 

The group $\Gamma$ can be extended to a group of isometries of $\hs$, for the details see \cite{Es.Ro}. The quotient $\hs / \Gamma$ is homeomorphic to $\re \times \si$, and its metric is given by 
$$dt^2 + \cosh^2(t)h_{\si}.$$ 
 
The surfaces $\{t\}\times \si$ are totally umbilical with constant mean curvature ${\rm tgh} (t)$, with respect to the unit normal $-\partial_t$. Fix $H \in (0,1)$, and define
\begin{equation}
t_H = {\rm arctgh} \ H.
\end{equation}
Applying a diffeomorphism we can rewrite the metric in $\hs / \Gamma $ as
\begin{equation}\label{Hmetric}
g_{H} = dt^2 + \cosh ^2 \left(t_H - t\right)h_{\si}.
\end{equation}
Now, $\{0\}\times \si$ has constant mean curvature equal to $H$.

\begin{teo}[Espinar-Rosenberg, \cite{Es.Ro}]
Let $(M^3,g)$ be a complete oriented Riemannian $3$-manifold whose sectional curvatures satisfy $K_{M} \leq -1$. Let $\si \subset M$ be an embedded closed oriented $H$-surface, $0 < H < 1$. Then
\begin{equation}\label{ER}
|\si| \leq \displaystyle\frac{2\pi|\chi(\si)|}{1 - H^{2}}
\end{equation}
Moreover if the equality holds and $\si$ separates the ambient space, then there are a neighbourhood $\mathcal{U}$ of $\si$ in $(M,g)$, and an isometry 
$$\Phi:\bigl([0,t_H)\times\si,g_{H}\bigr) \to (\mathcal{U},g),$$ 
such that $\Phi(0,\si) = \si$, where $g_{H}$ is the metric \eqref{Hmetric}.  Moreover, $\Phi$ can be extended to $\{t_H\}\times\si$, in such a way that $\si_m := \Phi(t_H,\si)$ is a totally geodesic surface, and $\si$ is a covering space of $\si_m$.
\end{teo}

It is important to remark that the area estimate \eqref{ER} also holds if $\si$ is a minimal surface, but the proof of the rigidity on the case of equality does not work if $H = 0$. On the last section we will show by a counter-example that the rigidity is not valid on this case. 

In \cite{Es.Ro}, Espinar and Rosenberg also obtained versions of the previous two theorems on the case of compact {\it capilarity} $H$-{\it surfaces} on Riemannian $3$-manifolds with boundary. 
\\

Observe that all the theorems stated concern compact surfaces. We are interested on the case where one have a complete non compact $H$-surface of finite area inside a $3$-manifold. As far as we know there is no area-rigidity theorem in this situation. We can now state our main results.
\begin{prop}\label{riginf}
Let $(M,g)$ be a complete oriented $3$-manifold with sectional curvature $K_{M} \leq -1$. Let $\si \subset M$ be a complete non-compact oriented $H$-surface, $0 \leq H < 1$, with finite topology and finite area. Then:
\begin{enumerate}
\item $\si$ has finite total curvature and \begin{equation}\label{GBF}
\int_{\si}K_{\si} \ d\mu = 2\pi\chi(\si);
\end{equation}
\item The area of $\si$ is estimated as \begin{equation}\label{area.est}
|\si| \leq \frac{2\pi|\chi(\si)|}{1 - H^{2}}.
\end{equation}
\end{enumerate}
Moreover, the equality on \eqref{area.est} is valid if, and only if, $\si$ is totally umbilic, $K_{M} \equiv -1$ along $\si$, and $K_{\si} \equiv H^{2} - 1$.
\end{prop}

We highlight that some versions of proposition \eqref{riginf} were already proved in \cite{C.H.R} and \cite{Me.Ra}.
\\

Suppose $H \in (0,1)$. In the case of equality in equation \eqref{area.est}, and under additional assumptions on $\si$ and $(M,g)$, a neighbourhood of the mean convex side of $\si$ must be locally isometric to a Fuchsian manifold. As we shall see in section \ref{couexa}, if $\si$ is a minimal surface one does not have rigidity in $(M,g)$. Precisely, we have:

\begin{thm}\label{loc.rig}
Let $(M^3,g)$ be a complete oriented Riemannian $3$-manifold whose sectional curvatures satisfy
$-a^{2} \leq K_{M} \leq -1$. Let $\si \subset M$ be an embedded complete non-compact oriented $H$-surface , $0 < H < 1$, with finite topology. Suppose that
$$|\si| = \displaystyle\frac{2\pi|\chi(\si)|}{1 - H^{2}}.$$ Then there are a neighbourhood $\mathcal{U}$ of $\si$ in $(M,g)$, and a local isometry 
$$\Phi:\bigl([0,t_H)\times\si,g_{H}\bigr) \to (\mathcal{U},g),$$ 
such that $\Phi(0,\si) = \si$, where $g_{H}$ is the metric \eqref{Hmetric}.
\end{thm}

\medskip

An important remark about the previous theorem is that the map $\Phi$ is is not necessarily a diffeomorphism. When $\si$ is closed and embedded in a $3$-manifold, it has an embedded tubular neighbourhood consisting of normal geodesics of same length, but in the non-compact case this is not always true even if one assumes that the immersion is proper. Consider for example the following situation:
\\\\
{\bf Example 1}:\label{bad} Let us review a construction done in \cite{Me.Ra}. Let $\mathcal{M}$ be a complete hyperbolic $3$-manifold of finite volume that contains a properly embedded totally geodesic thrice-punctured $2$-sphere $S$, the existence of such a manifold is proved in \cite{A}.
 
Let $N$ be a unit normal vector to $S$. Then, for any $t > 0$, the equidistant $S(t)=\{\exp_{x}\bigl(tN(x)\bigr); x\in S\}$ is the image
of a proper immersion of $S$ into $\mathcal{M}$
of constant mean curvature. Moreover, there is some $\varepsilon_{\mathcal{M}} > 0$ such that, for $t < \varepsilon_{\mathcal{M}}$, $S(t)$ is also properly embedded.

However, $S(t)$ intersects $S$ for all $t$, hence the map $\Psi: [0,\varepsilon_{\mathcal{M}})\times S \to \mathcal{M},$ $$\Psi(t,x) = \exp_{x}\bigl( (\varepsilon_{\mathcal{M}} - t) N(x)\bigr)$$ 
is a local diffeomorphism, but is not injective.

It is worth mention that a more general construction allows the totally geodesic surface to have other topology, \cite{A.M.R}.
\\\\

In view of this problem we will add more assumptions on $\si$ which assure that the map $\Phi$ on theorem \ref{loc.rig} is in fact a diffeomorphism, and so is an isometry. For state this result let us fix some notation. Consider a surface $\si$ properly embedded in a complete Riemannian $3$-manifold $(M,g)$ and denote by $N$ a unit normal along $\si$. Following chapter $8$ of \cite{G}, define the {\it minimal focal distance} of $\si$ in $M$ as
$$m_{f}(\si,N) := \inf_{x \in \si}\sup\bigl\{t > 0; \dist\bigl(\exp_{x}(tN(x)),\si\bigr) = t\bigr\}.$$
If $m_{f}(\si,N) > 0$, then every geodesic of the form $\gamma_{x}(t) = \exp_{x}(tN(x))$, $x \in \si$, is minimizing on the interval $[0,m_{f}(\si,N)]$, therefore $\si$ has an embedded tubular neighbourhood in $M$ consisting of normal geodesics of same length. 

Observe that the definition of $m_{f}(\si,N)$ depends on the choice of a unit normal. In the case of theorem \ref{loc.rig} we have $H_{\si} > 0$, hence the mean curvature vector $\vec{H}_{\si}$ has a well defined direction. So in this case we will choose $N$ pointing into the direction of $\vec{H}_{\si}$ and denote $m_{f}(\si) = m_{f}(\si,N)$. Following this we have:
\medskip

\begin{cor}\label{iso.rig}
Let $(M^3,g)$ and $\si$ satisfy the hypothesis of theorem \ref{loc.rig}. Suppose that $\si$ separates and $m_{f}(\si) > 0$. Then the map $\Phi$ of theorem \ref{loc.rig} is an diffeomorphism, and hence an isometry. Moreover, $\Phi$ can be extended to $\{t_H\}\times\si$, in such a way that $\si_m := \Phi(t_H,\si)$ is a totally geodesic surface, and $\si$ is a covering space of $\si_m$.
\end{cor}
\noindent
{\bf Remark}: In general an hyperbolic metric on a surface of finite topology is not unique. So in theorem \ref{loc.rig} and corollary \ref{iso.rig} we do not know what is the conformal structure of $\si$. However, it is well known that a thrice-punctured sphere has a unique hyperbolic metric up to diffeomorphism, so in this case we have only one possibility for the metric $h_{\si}$, and so for the metric $g_H$.
\\

Let us make some comments on the proof of the stated results. We follow the ideas of \cite{Ma.Ro} and \cite{Es.Ro}, however some difficulties arise due to the fact $\si$ is not compact. First, the analog of proposition \ref{riginf} in the case $\si$ is closed is proved using the {\it Gauss equation} and the {\it Gauss-Bonnet formula}, so we need to prove first the item \eqref{GBF} (which is interesting by itself) to obtain the area estimate \eqref{area.est}. On the proof of theorem \ref{loc.rig}, we work with the equidistants of $\si$. In the case $\si$ is closed, the geometry of these surfaces is uniformly controlled in a tubular neighbourhood of fixed radius, but in the non-compact case this is not necessarily true. In view of this, we added the hypothesis that the sectional curvature of $(M,g)$ is also bounded below, so we can use classical comparison theorems to control the geometry of the equidistants of $\si$. We also needed to obtain a version of the {\it first variation of area} for a family of non-compact surfaces with finite area, see equation \eqref{fvf}.
\\\\
{\bf Ackowledgements}: This work started while i worked at Instituto Nacional de Matem\'atica Pura e Aplicada (IMPA), under the fund {\it Programa de Capacita\c c\~ao Institucional PCI/MCTI-CNPq}. I would like to thank Jos\'e Espinar, Pedro Gaspar and Harold Rosenberg for discussions and comments on a early version of the manuscript. I also thank \'Alvaro Ramos, for discussions regarding example \eqref{bad}.

\section{Proof of the main results}

\begin{proof}[\bf Proof of Proposition \ref{riginf}]
Denote by $N_{\si}$ a unit vector field normal to $\si$ and by $A_{\si}$ the second fundamental form of $\si$ with respect to $N_{\si}$.

By the Gauss equation and the Arithmetic-Geometric Inequality, we obtain
\begin{equation}\label{eq1}
K _{\Sigma} = K_{M} + \det\bigl(A_{\si}\bigr) \leq -1 + H^{2} < 0.
\end{equation}

Since $\si$ has negative curvature, the curvature integrand $\int_{\si}K_{\si} \ d\sigma$ exists in $[-\infty,0)$. Also, $\si$ has finite area, then by theorem A in \cite{S} (see also theorem $12$ in \cite{H}), $\si$ has finite total curvature and it holds the 
\\\\
{\it Gauss-Bonnet Formula}:
\begin{equation}\label{eq2}
\int_{\si}K_{\si} \ d\mu = 2\pi\chi(\si).
\end{equation}

Integrating \eqref{eq1} and using \eqref{eq2} we have
$$(1 - H^{2})|\si| \leq - \int_{\si}K_{\si} \ d\mu = 2\pi|\chi(\si)|.$$

Moreover, equality holds if and only if, equality holds in \eqref{eq1}, that is, $K_\Sigma = H^2 - 1$, $K_{M} \equiv -1$ along $\Sigma$, and $\det\bigl(A_{\si}\bigr) \equiv H^{2}$ (which implies that $\si$ is totally umbilic).

\end{proof}

\begin{proof}[\bf Proof of theorem \ref{loc.rig}]
Let $N_{\si}$ be the unit normal along $\si$ pointing into the direction of the mean curvature vector $\vec{H}_{\si}$. Let $\Phi: [0,+\infty)\times\si \to M$ defined by $\Phi(t,x) = \exp_{x}\bigl(tN_{\si}(x)\bigr)$. Since $K_{M} \leq -1$ and $\si$ is totally umbilic with $H_{\si} < 1$, it follows from proposition $2.3$ in \cite{E} that $\si$ has no focal points, so the normal exponential map of $\si$ is a local diffeomorphism. Hence $\Phi$ is a local diffeomorphism. 

Using $\Phi$, we pull back the Riemannian metric of $M$ to $[0,+\infty)\times\si$. By the Gauss lemma, the vector field $\partial_{t}$ is everywhere normal to the surface $\si_{t} := \{t\}\times\si$, so this new metric can be written as
$$\Phi^{*}g = dt^{2} + \sigma_{t},$$
where $\{\sigma_{t}\}$ is a smooth family of metrics on $\si$.  

This construction makes of $\Phi$ a local isometry, so $\si_{0}$ is a totally umbilic $H$-surface on the metric $\Phi^{*}g$ with area $\displaystyle\frac{2\pi|\chi(\si)|}{1 - H^{2}}$. Moreover $\si_{t}$ is an equidistant to $\si_{0}$.
\\

Fix the notation $t_H = {\rm arctgh} \, H $. We will prove that for any $t \in [0,t_H)$, we have $$\sigma_{t} = \cosh^{2}(t_H - t)h_{\si},$$ where $h_{\si}$ is a hyperbolic metric of curvature $-1$, such that, $$g|_{\si} = \cosh^{2}(t_H)h_{\si}.$$

Denote by $A_{t}$ and $H_{t}$ respectively the second fundamental form and the mean curvature function of $\si_{t}$ with respect to the unit normal $\partial_{t}$. We also define $\lambda_{t}(x) \geq 0$, such that $H_{t} + \lambda_{t}$ and $H_{t} - \lambda_{t}$ are the principal curvatures of $\si_t$ at $(x,t)$. We notice that $\lambda_{t}(x) = 0$ if, and only if, $\si_t$ is umbilical at the point $(x,t)$.
\\

Consider the following pairs\\ 
{\bf Model 1}: $\bigl(\re\times\si,dt^2 + \cosh ^2 \left(t_H - t\right)h_{\si}\bigr)$, $\tilde\si_1 = \{0\}\times\si$\label{Mod1}.\\\\
{\bf Model 2}: $\bigl(\re\times\si,dt^2 + a^{2}\cosh ^2 \left(t_{a} - at\right)h_{\si}\bigr)$, $\tilde\si_2 = \{0\}\times\si$\label{Mod2}, where $t_{a} = {\rm arctgh} \ (a^{-1}H)$.\\

The model $1$ is just the Fuchsian model of curvature $-1$ we construct in the introduction and the equidistants of $\tilde\si_1$ have second fundamental form equal to 
$${\rm tgh}\left(t_H - t\right)\cdot h_{\si}.$$
Observe that for $t = 0$ we have the second fundamental form of $\tilde\si_1$ which is equal to $H\cdot h_{\si}.$
 
Furthermore, the model $2$ is a variation of this where the Riemannian manifold has constant sectional curvature $-a^2$ and the equidistants of $\tilde\si_2$ have second fundamental form equal to 
$$a\cdot{\rm tgh}\left(t_a - at\right)\cdot h_{\si}.$$
Observe that for $t = 0$ we have the second fundamental form of $\tilde\si_2$ which is equal to $H\cdot h_{\si}.$ 

Since $\si_{0}$ is totally umbilic and
\begin{equation}\label{eq3}
-a^{2} \leq K_{M} \leq -1,
\end{equation}
using the comparison theorem $3.1$ in \cite{E} (we compare $(M,g)$ with Model 1 and Model 2) we obtain 
\begin{equation}\label{eq4}
{\rm tgh}\left(t_H - t\right)\cdot h_{\si} \leq A_{t} \leq a\cdot{\rm tgh}\left(t_a - at\right)\cdot h_{\si},
\end{equation} 
where the inequalities are in the sense of quadratic forms.
So the functions $H_{t}$, $\lambda_{t}$ are uniformly bounded independently of $t$. Moreover, the theorem $3.2$ of \cite{He.Ka} implies that
\begin{equation}\label{eq5}
|\si_{t}| \leq \cosh(t_a - at)|\si|,
\end{equation}
so, if $t \in [0,\epsilon]$, the areas $|\si_{t}|$ are uniformly bounded independently of $t$.

Finally, we have the Gauss equation 
\begin{equation}\label{eq6}
K_{\si_t} = K_{M}(t) + (H_{t} + \lambda_{t})(H_{t} - \lambda_{t}),
\end{equation}
which together with the inequalities \eqref{eq3} and \eqref{eq4} implies that $K_{\si_t}$ is bounded. In particular the negative part $K_{\si_t}^{-}$ is bounded, so using that $\si_t$ has finite area we conclude that $\displaystyle\int_{\si_t}K_{\si_t}^{-}d\mu_t$ is finite. Then, by theorem $1$ in \cite{W1} (see also \cite{H}), it follows that $\si_t$ has finite total curvature. Using again theorem A in \cite{S}, we also obtain 
\begin{equation}\label{eq7}
\int_{\si_t} K_{\si_t} \ d\mu_{t} = 2\pi\chi(\si).
\end{equation}

Observe that the equation \eqref{eq6} implies
\begin{equation}\label{eq8}
K_{\si_t} \leq -1 + H_{t}^{2} - \lambda_{t}^{2}.
\end{equation}
Since $\si_t$ has finite area, and the functions $H_{t}^{2}$ and $\lambda_{t}^{2}$ are bounded and non-negative, its respective integrals over $\si_t$ exist and are finite. Integrating inequality \eqref{eq8} and using equation \eqref{eq7}, we have:
\begin{equation}\label{eq9}
2\pi\chi(\si) = \int_{\si_t} K_{\si_t} \ d\mu_{t} \leq - |\si_{t}| + \int_{\Sigma_t} H_{t}^{2}\ d\mu_{t} - \int_{\Sigma_t} \lambda_{t}^{2}\ d\mu_{t}.
\end{equation}

Thus we obtain the inequality
\begin{equation}\label{eq10}
\int_{\Sigma_t} \lambda_{t}^{2} \ d\mu_{t}  \leq - |\si_{t}| - 2\pi\chi(\si) + \int_{\Sigma_t} H_{t}^{2} \ d\mu_{t}.
\end{equation}

In the following, we denote by $F(t)$ the right hand side of inequality \eqref{eq10}.
\\\\
{\bf Claim}: F vanishes on $[0,t_H]$.
\begin{proof}[Proof of Claim 1]

First note that since $\si_0$ is a $H$-surface with area $\displaystyle\frac{2\pi|\chi(\si)|}{1 - H^{2}}$, we have $F(0) = 0$.

Now, observe that $\{\si_t; t \in \re\}$ is a normal variation of $\si_0$ with variational vector field $X = \partial_{t}$. Let $\{B_{n}; n \in \mathbb{N}\}$ be a exhaustion of $\si$ by compact domains with smooth boundary, and denote by $B_{n,t}$ the image of $B_{n}$ by the flow of $\partial_{t}$ which is inside $\si_t$. By the first variation formula
$$
\frac{d}{dt}|B_{n,t}| = - 2\int_{B_{n,t}} H_{t} \ d\mu_{t},
$$
observe that there is no boundary term since the variational vector field $X$ is everywhere normal to $\si_t$. 

By inequality \eqref{eq4}
\begin{equation}\label{eq11}
0 \leq H_t \leq a, \ \ \forall t \in [0,t_H],
\end{equation}
which together with the fact $\si_t$ has finite area implies that $\int_{\si_t}H_{t}\ d\mu_{t}$ exist and is finite.

Furthermore, using equations \eqref{eq4} and \eqref{eq5}, we see there are constants $c,\tilde{c} > 0$ independent of $t$ such that
$$
\biggl|\int_{\si_{t}} H_{t} \ d\mu_{t} - \int_{B_{n,t}} H_{t} \ d\mu_{t}\biggr| = \biggl|\int_{\si_{t}\backslash B_{n,t}} H_{t} \ d\mu_{t}\biggr| \leq c |\si_t\backslash B_{n,t}| \leq \tilde{c}|\si\backslash B_n|.
$$

The last term on the above formula converges to $0$ when $n$ goes to $+\infty$, hence $\bigl\{\frac{d}{dt}|B_{n,t}|\bigr\}_{n \in \mathbb{N}}$ converges uniformly to $-2\int_{\si_t} H_{t} \ d\mu_{t}$. Moreover, $\bigl\{|B_{n,t}|\bigr\}_{n \in \mathbb{N}}$ converges to $|\si_{t}|$, for any $t$. By the criterion of differentiability for limits of sequences of functions, we conclude that $|\si_{t}|$ is differentiable for $t \in [0,t_H]$ and we obtain the
\\\\
{\it First Variation of Area}: 
\begin{equation}\label{fvf}
\frac{d}{dt}|\si_{t}| = - 2\int_{\si_{t}} H_{t} \ d\mu_{t}.
\end{equation}

Now, define $F_n(t) = - |B_{n,t}| - 2\pi\chi(\si) + \displaystyle\int_{B_{n,t}} H_{t}^{2} \ d\mu_{t}.$ By the first variation of mean curvature 
\begin{equation}\label{fvmc}
\frac{\partial H_t}{\partial t} = \frac{1}{2}\bigl(\ric_{M}(\partial_{t},\partial_{t}) + |A_{t}|^{2}\bigr) = \frac{1}{2}\ric_{M}(\partial_{t},\partial_{t}) + H_{t}^{2} + \lambda_{t}^{2}.
\end{equation} 
Thus,
\begin{eqnarray*}
\frac{d F_n}{dt} &=& \int _{B_{n,t}} 2 H_t \ d\mu_{t} + \int _{B_{n,t}} \biggl(2 H_t \frac{\partial H_t}{\partial t} - 2H_{t}^{3}\biggr)d\mu_{t}\\\\
&=& \int _{B_{n,t}} H_t\bigl(2 + \ric_{M}(\partial_{t},\partial_{t})\bigr)d\mu_{t} + 2\int _{B_{n,t}} H_t\lambda_{t}^{2}\ d\mu_{t},
\end{eqnarray*}
for all $t \in [0,t_H]$. 

We obtain from inequalities \eqref{eq3} and \eqref{eq11} that 
\begin{equation}\label{eq12}
2a(1 - a^{2}) \leq H_t\bigl(2 + \ric_{M}(\partial_{t},\partial_{t})\bigr) \leq 0, \ \ \forall t \in [0,t_H],
\end{equation}
and
\begin{equation}\label{eq13}
0 \leq H_t\lambda_{t}^{2} \leq C, \ \ \forall t \in [0,t_H],
\end{equation}
for some constant $C > 0$.
Therefore, the integrals $\displaystyle\int _{\si_t} H_t\bigl(2 + \ric_{M}(\partial_{t},\partial_{t})\bigr)\ d\mu_{t}$ and $\displaystyle\int _{\si_t} H_t\lambda_{t}^{2}\ d\mu_{t}$ exist and are finite. Arguing as we did to prove the {\it First Variation of Area} we conclude that the function $F$ is differentiable for $t \in [0,t_H]$ and
\begin{equation}\label{eq14}
\frac{d F}{dt} = \int _{\si_t} H_t\bigl(2 + \ric_{M}(\partial_{t},\partial_{t})\bigr)d\mu_{t} + 2\int_{\si_t} H_t\lambda_{t}^{2}\ d\mu_{t} \leq  2C\int_{\si_t} \lambda_{t}^{2}\ d\mu_{t},
\end{equation} 
where in the last part we used inequalities \eqref{eq12} and \eqref{eq13}.

Since $$\int _{\si_t} \lambda_{t}^{2} \ d\mu_{t} \leq F(t),$$
we obtain $\displaystyle\frac{d F}{dt} \leq 2CF$ in $[0,t_H]$.

By Gronwall's Lemma, we get $F(t) \leq F(0)e^{2Ct}$ on $[0,t_H]$. But $F(0) = 0$, so $F(t) \leq 0$ on $[0,t_H]$. However, by inequality \eqref{eq4} we have $F(t) \geq 0, \forall t$. Then it follows that $F \equiv 0$.
\end{proof}

Since $F \equiv 0$ all the inequalities in formula \eqref{eq14} are in fact equalities, so we obtain:
\begin{itemize}
\item [(a)] $\Sigma _t$ is umbilic;
\item [(b)] $\ric_{M}(\partial_{t},\partial_{t}) + 2 \equiv 0$, which implies that the sectional curvature of any $2$-plane orthogonal to $\Sigma_t$ is equals to $-1$;
\item [(c)] $\displaystyle\frac{\partial H_t}{\partial t} = H_t^2 -1 $.
\end{itemize}

On the one hand, by items (a) and (c) we have 
\begin{equation}\label{HNeg}
H_t (x) = {\rm tgh}(t_H - t) \text{ for all } x \in \si ,
\end{equation}
that is, $\si_t$ has constant mean curvature for all $t \in [0,t_H]$.

On the other hand, by the {\it First Variation of the Area} \eqref{fvf}, we have 
$$\frac{d}{dt} \abs{\Sigma _t} = 2 \, {\rm tgh}\left(t_H - t\right) \abs{\Sigma _t} \text{ for all } t \in [0,t_H), $$
which implies that  
\begin{equation}\label{AreaNeg}
 \abs{\Sigma _t} =  2\pi \abs{\chi(\Sigma)}\cosh^2 \left(t_H - t\right) \text{ for all } t \in [0,t_H].
\end{equation}

It follows from inequality \eqref{eq9} and from equations \eqref{HNeg} and \eqref{AreaNeg} that 
\begin{equation*}
2\pi \chi(\Sigma) = \int _{\Sigma _t} K_{\Sigma _t}  \leq (H_t^2 -1) \abs{\Sigma _t} = - 2 \pi \abs{\chi(\Sigma)},
\end{equation*}
hence the sectional curvature of the $2$-plane tangent to $\Sigma _t$ is equal to $-1$ for all $t \in [0,t_H]$. In view of all this we conclude that
$$K_{\Sigma _t} = {\rm tgh}^2 (t_H - t) - 1= -\cosh^{-2} \left(t_H - t\right), \ \forall t \in [0,t_H].$$

Therefore the metric $dt^2 + \sigma_t$ is hyperbolic. Since, $\si$ is totally umbilical with constant mean curvature $H$, we conclude that
$$\sigma_{t} = \cosh^{2}(t_H - t)h_{\si}, \ \forall t \in [0,t_H].$$

\end{proof}

\begin{proof}[\bf Proof of corolary \ref{iso.rig}]
If $\Phi$ is injective on $ [0, t_H] \times \si $ there is nothing to do, so suppose this is not true and define
 $$ \bar \epsilon = \sup\bigl\{\epsilon \in [0,t_H); \Phi \textrm{ is injective on } [0,\epsilon ) \times \si\bigr\}. $$

Since $m_{f}(\si) > 0$, we have $\bar \epsilon > 0$. Suppose that $\Phi$ is not injective on $ [0, t_H) \times \si $. Then there are distinct points $p,q \in \si $ such that $1$ or $2$ holds:
\begin{enumerate}
\item $\Phi(\bar \epsilon , p) = \Phi(0, p ) \text{ or } \Phi(0, q) $,
\item $\Phi(\bar \epsilon , p) = \Phi(\bar \epsilon , q) $.
\end{enumerate}

In case $1$, we will construct a closed curve $C$ meeting tranversely a equidistant $\si_{t_0}$ in exactly one point, where $0< t_0 < \bar{\epsilon}$, but this contradicts the fact that $\si$ separates $M$. To construct $C$ when $\Phi( \bar \epsilon , p) = \Phi (0, p)$, take $C(t) = \Phi(t, p)$, $0 \leq t  \leq  \bar \epsilon $. When $\Phi(\bar \epsilon , p) = \Phi (0, q)$, let $\Gamma$ be a curve joining $q$ to $p$ on $\si$, and define $C = \{\Gamma\}\cup\{\Phi(t,p)$; $0 \leq  t  \leq  \bar \epsilon\}$.  Thus case $1$ can not occur.

Consider case $2$, $\Phi( \bar \epsilon ,p) = \Phi(\bar \epsilon ,q) =: x_0$. Given $W \subset \si_0$ define $W_t = \{t\}\times W$. Since $\Phi$ is a local diffeomorphism, there exists neighbourhoods $U$ and $V$ of $p$ and $q$ (respectively) in $\si_0$, such that $[0, \bar \epsilon)\times U$ and $[0, \bar \epsilon)\times V$ are disjoint, $\Phi$ is injective on $[0, \bar \epsilon]\times U$ and on $[0, \bar \epsilon ]\times V$, and $x_0 \in \Phi\bigl(U_{\epsilon}\bigr)\cap\Phi\bigl(V_{\epsilon}\bigr)$.

There are two possibilities, depending on the direction of the mean curvature vectors, $\vec{H}_{\Phi(\bar \epsilon ,p)}$ and $\vec{H}_{\Phi(\bar \epsilon , q)}$. Suppose first that $\vec{H}_{\Phi(\bar \epsilon ,p)} = \vec{H}_{\Phi(\bar \epsilon , q)}$. Then $\Phi(U_{\bar\epsilon})$ is on one side of $\Phi\bigl(V_{\bar\epsilon}\bigr)$ at $x_0$ (where they are tangent) and they have the same mean curvature vector at this point, so by the maximum principle $\Phi(U_{\bar\epsilon}) = \Phi(V_{\bar\epsilon})$. Consequently, $\si$ is a non trivial covering space of the embedded surface $\Phi(\si_{\bar\epsilon})$. However, for $\tilde \epsilon$ near $\bar \epsilon $, $\Phi(\si_{\tilde\epsilon})$ is a graph over $\Phi(\si_{\bar\epsilon})$ in the trivial normal bundle of $\Phi(\si_{\bar\epsilon})$ in $M$. This is impossible when the covering space is non trivial.

Now suppose that $\vec{H}_{\Phi(\bar \epsilon ,p)} = - \vec{H}_{\Phi(\bar \epsilon , q)}$. At $x_0$, $\Phi(U_{\bar\epsilon})$ is strictly mean convex towards $\vec{H}_{\Phi(\bar \epsilon , p)}$  and $\Phi(\{\bar \epsilon\}\times V)$ is strictly mean convex towards $-\vec{H}_{\Phi(\bar \epsilon ,p)} $. Then $\Phi\bigl([0, \bar \epsilon)\times U\bigr)$ is on the mean concave side of $\Phi(\{\bar \epsilon\}\times U)$ near $x_0$, and $\Phi\bigl([0, \bar \epsilon)\times V\bigr)$ is on the mean concave side of $\Phi(V_{\bar\epsilon})$ near $x_0$. So, any small neighbourhood of $\Phi(\bar \epsilon , p)$ in $M$ contains an open set on the mean concave side of $\Phi\bigl(U_{\bar\epsilon}\bigr)$ which is contained in the mean concave side of $\Phi\bigl(V_{\bar\epsilon}\bigr) $. Hence $\Phi\bigl([0, \bar \epsilon)\times U\bigr)$ intersects $\Phi\bigl([0, \bar \epsilon)\times V\bigr)$ which contradicts the fact that $\Phi$ injective on $[0,\bar \epsilon) \times \si $.

Therefore $\Phi$ is injective on $ [0, t_H ) \times \si $. Thus $\Phi|_{[0, t_H ) \times \si }$ is a diffeomorphism. Moreover, we know that $\Phi(t_H,\cdot)$ is a immersion, and by the proof of theorem \ref{loc.rig}, the surface $\si_{t_H }$ is totally geodesic with respect to the metric $\Phi^{*}g$. Hence $\Phi\bigl(\si_{t_H}\bigr)$ is a totally geodesic surface. Arguing as in the case 2 discussed before we conclude that $\si$ is a covering space of $\Phi\bigl(\si_{t_H}\bigr)$.
\end{proof}

\section{Counterexample in the case $H = 0$}\label{couexa}

Let $S$ be a oriented surface which is either closed or non-compact with finite topology, endowed with a complete metric $h_{S}$ of constant curvature $-1$ and finite area. Consider the manifold $M = [0,+\infty)\times S$ endowed with the warped product metric $$g = dt^2 + (e^t - t)^{2}h_{S}.$$

Denote $f(t) = e^t - t$. Consider the surface $\si = \{0\}\times S$. Its induced metric is $g|_{\si} = h_{S}$ and its second fundamental form is $A_{\si} = -\frac{f'(0)}{f(0)}g \equiv 0$. So $\si$ is totally geodesic,
$$K_{\si} \equiv -1 \ \textrm{and} \ |\si| = 2\pi|\chi(S)|.$$

Denote $\xi = \partial_{t}$. Given a vector field $V$ on $M$, define the projection on $S$ as $V^S := V - g(V,\xi)\xi$. Following the convention in \cite{ON} to sign of the curvature, we have that the curvature tensor $R_{M}$ of M can be expressed in terms of the
curvature tensor $R_{S}$ of $S$ as (see \cite{ON}, prop.12.42)
\begin{eqnarray*}
R_{M}(X,Y)Z &=& R_{S}(X^S,Y^S)Z^S - \biggl(\frac{f'}{f}\biggr)^2\bigl[g(X,Z)Y - g(Y,Z)X\bigr]\\
&+& (\log f)''g(Z,\xi)\bigl[g(Y,\xi)X - g(X,\xi)Y\bigr]\\
&-& (\log f)''\bigl[g(X,Z)g(Y,\xi) - g(Y,Z)g(X,\xi)\bigr]\xi,
\end{eqnarray*}
where $X,Y,Z$ are smooth vector fields on $M$. Given a plane $\Pi \subset T_{(t,x)}M$, and an orthonormal basis $\{X,Y\}$ of $\Pi$, it follows that the sectional curvature of $\Pi$ is given by
\begin{eqnarray*}
K_{M}(\Pi) &=& \frac{K_{S} - (f')^2}{f^2} + \biggl[\frac{-K_{S} + (f')^2}{f^2} - \frac{f''}{f}\biggr]\bigl[g(X,\xi)^2 + g(Y,\xi)^2\bigr]\\
&=& \frac{-1 - (e^t - 1)^2}{(e^t - t)^2} + \biggl[\frac{1 + (e^t - 1)^2}{(e^t - t)^2} - \frac{e^t}{e^t - t}\biggr]\bigl[g(X,\xi)^2 + g(Y,\xi)^2\bigr]\\
&=& \frac{-1 - (e^t - 1)^2}{(e^t - t)^2} + \biggl[\frac{e^t(t - 2) + 2}{(e^t - t)^2}\biggr]\bigl[g(X,\xi)^2 + g(Y,\xi)^2\bigr].
\end{eqnarray*}

Since $\xi$ is unitary we have $0 \leq g(X,\xi)^2 + g(Y,\xi)^2 \leq 1$. Moreover, there is $0 < t_0 < 1$ such that the function $u$ defined by $u(t) = \displaystyle\frac{e^t(t - 2) + 2}{(e^t - t)^2}$ is non-positive, for all $t \in [0,t_0]$, and is positive for all $t > t_0$. So, if $t \in [0,t_0]$, the sectional curvature will be maximal when $g(X,\xi)^2 + g(Y,\xi)^2 = 0$, and in this case
$$K_{M}(\Pi) = \frac{-1 - (e^t - 1)^2}{(e^t - t)^2}.$$
Define $w(t) = \displaystyle\frac{-1 - (e^t - 1)^2}{(e^t - t)^2}$. We have
\begin{equation}\label{deriv}
w'(t) = \frac{2(e^t - t)(e^t - 1)}{(e^t - t)^4}\bigl[e^t(t - 2) + 2\bigr],
\end{equation}
so the sign of $w'(t)$ is determined by the sign of $e^t(t - 2) + 2$, hence $w'(t) \leq 0$ for $t \in [0,t_0]$, and thus $$K_{M}(\Pi) = w(t) \leq w(0) = -1.$$

On the other hand, if $t > t_0$, the sectional curvature will be maximal when $g(X,\xi)^2 + g(Y,\xi)^2 = 1$, and in this case
$$K_{M}(\Pi) = -\frac{e^t}{e^t - t} \leq -1,$$
where the inequality holds since $e^t \geq e^t - t, \forall t \geq 0$.

Now we want to prove that the sectional curvatures of $(M,g)$ are also bounded below. Observe that, if $t \in [0,t_0]$, the curvature will be minimal when $g(X,\xi)^2 + g(Y,\xi)^2 = 1$, i.e,
$$K_{M}(\Pi) = -\frac{e^t}{e^t - t},$$
since the left hand side defines a continuous function, we have that this function is bounded in $[0,t_0]$. If $t > t_0$, the curvature will be minimal when $g(X,\xi)^2 + g(Y,\xi)^2 = 0$, i.e,
$$K_{M}(\Pi) = \frac{-1 - (e^t - 1)^2}{(e^t - t)^2}.$$
It follows from equation \eqref{deriv} that $w$ is increasing on $(t_0,+\infty)$, and then $w(t) \geq w(t_0), \forall t > t_0$. 

Therefore $(M,g)$ satisfies $-a^2 \leq K_{M} \leq -1$, for some $a > 0$, but the metric $g$ is not hyperbolic. To obtain a complete manifold without boundary, just take the double of M, the metric $g$ in the two copies of $M$ will glue appropriately because $\si$ is totally geodesic. This construction provides an counter-example to theorem \eqref{loc.rig} and corollary \eqref{iso.rig} in the case $H = 0$.

On the other hand, consider the Riemannian metric $$\tilde{g} = dt^2 + l(t)^{2}h_{S},$$ on $\re\times S$, where
$$l(t) =  \left\{
\begin{array}{rl}
e^t - t & \text{if }, t \geq 0,\\\\
\cosh t & \text{if }, t < 0.
\end{array} \right..$$
Since the function $l$ is of class $C^2$ (but is not $C^3$), we obtain a $C^2$ metric which agrees with the model up to reach the totally geodesic slice, where then the metric agrees with our counter-example. 

Concerning this last construction, observe that the proofs of the main results only use $C^2$ properties of the metric. It is interesting to know if an example like this can exist with an $C^{\infty}$ metric.

$$$$

\noindent \textsc{Instituto de Matem\'atica e Estat\'istica, UERJ\\ Rua S\~ao Francisco Xavier, 524\\Pavilh\~ao Reitor Jo\~ao Lyra Filho, 6º andar - Bloco B\\ 20550-900, Rio de Janeiro-RJ, Brazil}
\\\\
\noindent {\it Email address}: \texttt{vanderson@ime.uerj.br}

\end{document}